\documentclass[11pt]{article}

\usepackage[utf8]{inputenc}
\usepackage{amsmath}
\usepackage[english]{babel}
\usepackage{amssymb}
\usepackage{amsthm}
\usepackage{mathtools}
\usepackage{esdiff}
\usepackage{microtype}
\usepackage{paralist}
\usepackage{skak}
\usepackage{bbm}
\usepackage[paper=a4paper,left=25mm,right=25mm,top=25mm,bottom=25mm]{geometry}
\usepackage{stmaryrd}
\usepackage{hyperref}
\usepackage{color} 
\usepackage{fancyhdr}
\usepackage{nomencl}
\usepackage{cite}
\usepackage{longtable}
\usepackage{thmtools}
\usepackage{cuted}
\usepackage{MnSymbol}
\declaretheoremstyle[headfont=\normalfont\bfseries]{bfthmstyle}

\declaretheoremstyle[headfont=\normalfont\bfseries,qed={$\diamond$}]{bfdefstyle}
\declaretheorem[sharenumber=Theorem,style=bfdefstyle]{Definition}
\declaretheorem[sharenumber=Theorem,style=bfdefstyle]{Definition-Proposition}
\declaretheorem[sharenumber=Theorem,style=bfdefstyle]{Definition-Lemma}

\newtheorem{definition-proposition}[Theorem]{Definition-Proposition} %avoid in the future
\newtheorem{definition-lemma}[Theorem]{Definition-Lemma} %avoid in the future

\declaretheoremstyle[headfont=\normalfont\bfseries,qed={$\diamond$}]{bfremstyle}

\declaretheoremstyle[
notefont=\normalfont, notebraces={}{},
headformat=\mathbb{N}UMBER~\mathbb{N}AME~\mathbb{N}OTE
]{nopar}

 \numberwithin{equation}{section}

\makenomenclature

\author{Bernhard Maschke, Jonas Kirchhoff}
\title{Port maps of Irreversible Port Hamiltonian Systems}
\date{\today}

\begin{document}
\setlength{\parindent}{0em}
\pagestyle{fancy}
\lhead{Bernhard Maschke, Jonas Kirchhoff}
\rhead{Port Maps}

\maketitle

\paragraph{Abstract} Irreversible Port Hamiltonian Systems are departure of Port Hamiltonian
Systems as they are generated not only by a Hamiltonian function but
also by an entropy function and defined with respect to a quasi-Poisson
bracket which embeds the definition of the irreversible phenomena
taking place in the system. However the port map, consisting in the
input map and the output map were poorly justified and lacked any
physical consistency. In this paper, we suggest a novel definition
of the port maps which allows to recover not only the energy balance
equation (when the Hamiltonian equals the total energy of the system)
but also a entropy balance equation including the irreversible entropy
creation at the interface (the port) of the system in addition to
the entropy creation term due to internal irreversible phenomena.\vspace{-1mm}
 
\paragraph{Keywords} Port Hamiltonian Systems, Nonlinear Systems, Irreversible Thermodynamics, Energy and Entropy based Modelling, Geometrical Methods.\vspace{-1mm}

\vfill
\par\noindent\rule{5cm}{0.4pt}\\
\begin{footnotesize}
Corresponding author: Bernhard Maschke\\[1em]
Bernhard Maschke\\
Univ. Lyon, Université Claude Bernard Lyon 1, CNRS, LAGEPP UMR 5007, France\\
E-mail: bernhard.maschke@univ-lyon1.fr\\[1em]
Jonas Kirchhoff\\
Institut für Mathematik, Technische Universität Ilmenau, Weimarer Stra\ss e 25, 98693 Ilmenau, Germany\\
E-mail: jonas.kirchhoff@tu-ilmenau.de

\end{footnotesize}

\newpage

\section{Introduction}

Irreversible Port Hamiltonian Systems are a class of Port Hamiltonian
Systems tailored to represent physical systems subject to dissipative
phenomena and including entropy balance equations \cite{Ramirez_ChemEngSci13}.
Their dynamics is generated by a Hamiltonian function which, for physical
systems, is equal to the total energy and defined with respect to
a quasi-Poisson bracket. In this sense it depart from most other geometrical
formulations of irreversible thermodynamic systems where the irreversible
phenomena are formulated as a (pseudo-)gradient system such as in
\cite{GrmelaOttinger97,Morrison86}. This quasi-Poisson bracket has
a very particular structure as it is defined as a Poisson bracket,
premultiplied by a product of functions depending not only on the
state variable but also on the differential of the Hamiltonian. This
bracket is the object of a companion paper \cite{Kirchhoff_submIFAC_WC22}.

An extension to systems encompassing as well irreversible as reversible
coupling phenomena has been suggested in \cite{Ramirez_EJC13} and
considered in numerous examples, and recent extension to distributed
parameter systems has been given \cite{Ramirez_ChemEngSci_2022,Ramirez_Entropy22}.
Taking explicitely account of the entropy balance equation in the
Irreversible Port Hamiltonian formulation, has given rise to novel
nonlinear controller based on shaping the closed-loop entropy production
also to a novel optimal control synthesis \cite{Maschke_ArXiv22_optContrIPHS}. 

However it is remarkable that in the various publications on Irreversible
Port Hamiltonian Systems, one finds various definitions of the input
map and sometimes no definition of the conjugated output variable.
In this paper, we shall be interested in the case when the input map
is affine in the control variables, as defined in \cite{Ramirez_ChemEngSci13},
which corresponds to assume that the interface of the system with
its environment, is subject to an irreversible phenomenon such as
heat conduction for instance. And we shall show how to derive port
maps which are compatible with the quasi-Poisson bracket of Irreversible
Port Hamiltonian Systems and illustrate the constructionon a very
elementary example.

In the section 2, we shall recall the definition of
Irreversible Port Hamiltonian Systems and state precisely the motivation
of the paper. In the section 3, we
shall suggest a more precise definition of the affine input-map associated
with irreversible (or dissipative) interfaces and derive this definition
for the general quasi-Poisson brackets of IPHS.

\section{Irreversible Port Hamiltonian Systems and problem
statement}\label{sec:IPHS}

\subsection{Reminder on Irreversible Port Hamiltonian Systems}

Let us first, briefly recall the definition of Irreversible Port Hamiltonian
System \cite{Ramirez_ChemEngSci13} . 
\begin{Definition}
\label{def:IPHS} \cite{Ramirez_ChemEngSci13} An \emph{Irreversible
Port Hamiltonian Systems} (IPHS) is the nonlinear control system 
\begin{equation}
\frac{dx}{dt}=\gamma\left(x,\tfrac{\partial U}{\partial x}\right)\left\{ S,U\right\} _{J}J\frac{\partial H}{\partial x}(x)+W\left(x,\tfrac{\partial H}{\partial x}\right)+g\left(x,\tfrac{\partial H}{\partial x}\right)u,\label{eq:IPHS}
\end{equation}
where $x\left(t\right)\in\mathbb{R}^{n}$ is the state vector, $u\left(t\right)\in\mathbb{R}^{m}$
is the control input, and defined by 

(i) two (smooth) real functions called \emph{Hamiltonian function}
$H(x)\in C^{\infty}(\mathbb{R}^{n})$ and \emph{entropy function}
$S(x)\in C^{\infty}(\mathbb{R}^{n})$, 

(ii) the \emph{structure matrix} $J\in\mathbb{R}^{n}\times\mathbb{R}^{n}$
which is constant and skew-symmetric defining a \emph{Poisson bracket}
\emph{\cite{Libermann_marle87}}
\begin{equation}
\left\{ S,H\right\} _{J}=\frac{\partial S}{\partial x}^{\top}\left(x\right)J\frac{\partial H}{\partial x}\left(x\right)\label{eq:PoissonBracket}
\end{equation}

(iii) a real function $\gamma(x,\frac{\partial U}{\partial x})=\hat{\gamma}(x):C^{\infty}(\mathbb{R}^{n})$,
strictly positive function of the states and co-states

(iv) the vector field $W(x,\frac{\partial U}{\partial x})\in\mathbb{R}^{n}$
and matrix field $g(x,\frac{\partial U}{\partial x})\in\mathbb{R}^{n\times m}$
associated with the input map. 
\end{Definition}

Firstly by the skew-symmetry of $J$, it follows that the energy is
a conserved quantity which obeys the following balance equation
\begin{equation}
\begin{split}\frac{dH}{dt} & =\frac{\partial H}{\partial x}^{\top}\left(\gamma\left(x,\tfrac{\partial U}{\partial x}\right)\left\{ S,H\right\} _{J}J\frac{\partial H}{\partial x}\right)+\frac{\partial H}{\partial x}^{\top}\left(W+gu\right)\\
 & =\gamma\left(x,\tfrac{\partial U}{\partial x}\right)\left\{ S,H\right\} _{J}\underbrace{\left\{ H,H\right\} _{J}}_{=0}+\frac{\partial H}{\partial x}^{\top}\left(W+gu\right)\\
 & =\frac{\partial H}{\partial x}^{\top}\left(W+gu\right)
\end{split}
\label{eq:EnergyBalanceIPHS}
\end{equation}
 It may be observed that the right-hand side of the energy balance
equation (\ref{eq:EnergyBalanceIPHS}) which is the power incoming
the system, may not be interpred as the product of the input and a
conjuguated port-output variable, due to the affine nature of the
input map. Hence the system may not be qualified as being impedance
passive as it is the case for Port Hamiltonian Systems.

The entropy balance of the system is given by 
\begin{equation}
\begin{split}\frac{dS}{dt} & =\frac{\partial S}{\partial x}^{\top}\left(\gamma\left(x,\tfrac{\partial U}{\partial x}\right)\left\{ S,H\right\} _{J}J\frac{\partial H}{\partial x}\right)+\frac{\partial H}{\partial x}^{\top}\left(W+gu\right)\\
 & =\underbrace{\gamma\left(x,\tfrac{\partial U}{\partial x}\right)\left\{ S,H\right\} _{J}^{2}}_{=\sigma_{int}}+\frac{\partial S}{\partial x}^{\top}\left(W+gu\right)
\end{split}
\label{eq:EntropyBalanceIPHS}
\end{equation}
where $\sigma_{int}$ is called the \textit{internal entropy production}
and is positive and, due to the strict positivity of the function
$\gamma$, zero if and only if $\left\{ S,H\right\} _{J}=0$. For
numerous physical systems $\left\{ S,H\right\} _{J}$ appears to be
the driving force of the irreversible phenomena which implies that
if it is zero, the system is at equilibrium and no irreversible phenomenon
takes place: the entropy production is hence also zero. The structure
of the drift dynamics, i.e. the two functions $\gamma$ and $\left\{ S,H\right\} _{J}$
premultiplying the structure matrix $J$, ensures that the second
principle of the Thermodynamics is satisfied: the drift dynamics yieds
an increase of the entropy of the system equal to internal entropy
production. But the second term of the entropy balance equation (\ref{eq:EntropyBalanceIPHS}),
the entropy flow coming from the environment, again may not be interpreted
as a impedance type supply. Moreover, the completely free choice of
$W$ and $g$ do not reflect that the interface has also to obey the
two principles of Thermodynamics.

In order to illustrate the latter observation, let us consider the
following elementary example.

\subsection{Heat conduction in 2 compartement system \cite{Ramirez_ChemEngSci13} }

Consider two compartment, indexed by 1 and 2 (constaining for instance
two ideal gases), which may interact through a conducting wall and
the compartment 2 interacting through a heat conducting wall with
a controlled thermostat at temperature $T_{e}(t)$ . The dynamics
of this system is given by the entropy balance equations of each compartment
\[
{\small\frac{d}{dt}\begin{bmatrix}S_{1}\\
S_{2}
\end{bmatrix}=\lambda\begin{bmatrix}\frac{T_{2}(S_{2})}{T_{1}(S_{1})}-1\\
\frac{T_{1}(S_{1})}{T_{2}(S_{2})}-1
\end{bmatrix}+\lambda_{e}\begin{bmatrix}0\\
\frac{T_{e}(t)}{T_{2}(S_{2})}-1
\end{bmatrix}}
\]
where $S_{1}$ and $S_{2}$ are the entropies of compartement 1 and
2, $T_{e}(t)>0$ is the controlled thermostat temperature and $\lambda>0$
and $\lambda_{e}>0$ denotes Fourier's heat conduction coefficients
of the two walls. Assuming that the two compartments contain a pure
ideal gas and that they undergo no deformation, and are closed, the
temperatures may be modelled as exponential functions of the entropies
$T(S_{i})=T_{0}\exp\left(\tfrac{S_{i}}{c_{i}}\right)$, where $T_{0}$
and $c_{i}$ are constants.

This system may be written as a IPHS with state vector $x^{\top}=\left(\begin{array}{cc}
S_{1} & S_{2}\end{array}\right)$ being the entropy of each compartement, Hamiltonian function being
the total internal energy $H\left(x_{1},x_{2}\right)=U_{1}(x_{1})+U_{2}(x_{2})$
such that $T_{i}\left(x_{i}\right)=\frac{\partial H}{\partial x_{i}}>0$,
entropy function being the total entropy $S\left(x_{1},x_{2}\right)=x_{1}+x_{2}$.
The structure matrix is the symplectic matrix $J=\left[\begin{smallmatrix}0 & 1\\
-1 & 0
\end{smallmatrix}\right]$ such that 
\[
\{S,U\}_{J}=\frac{\partial S}{\partial x}^{\top}J\frac{\partial U}{\partial x}=\begin{bmatrix}1\\
1
\end{bmatrix}^{\top}\begin{bmatrix}0 & 1\\
-1 & 0
\end{bmatrix}\begin{bmatrix}T_{1}\\
T_{2}
\end{bmatrix}=T_{2}-T_{1}
\]
 is indeed the driving force of the heat conduction. And the Fourier's
law is contained in the definition of the positive function $\gamma\left(x\right)=\frac{\lambda}{T_{1}T_{2}}$.
The input variable is $u=T_{e}$ and the input map defined by the
vector $g=\left(\begin{array}{c}
0\\
\frac{\lambda_{e}}{T_{2}(S_{2})}
\end{array}\right)$ and $W=\left(\begin{array}{c}
0\\
-\lambda_{e}
\end{array}\right)$.

Computing the total entropy balance equation (\ref{eq:EntropyBalanceIPHS}),
one obtains 
\[
\frac{dS}{dt}=\underbrace{\lambda\frac{\left(T_{2}-T_{2}\right)^{2}}{T_{1}T_{2}}}_{=\sigma_{int}}+\frac{\lambda_{e}}{T_{2}}\left(u-T_{2}\right)
\]
 where the entropy creation $\sigma_{e}$ due to the heat conduction
at the external wall is not apparent whereas by substracting and adding
the thermostat entropy flow entering the system $f_{port}^{S}=\frac{\left(u-T_{2}\right)}{u}$,
one obtains an expression where it appears
\begin{equation}
\frac{dS}{dt}=\frac{\lambda}{T_{1}T_{2}}\left(T_{2}-T_{2}\right)^{2}+\underbrace{\frac{\lambda_{e}}{uT_{2}}\left(u-T_{2}\right)^{2}}_{=\sigma_{port}\geq0}+f_{port}^{S}\label{eq:IrrevesibelEntrCrea_WallExt}
\end{equation}

\subsection{Problem statement}

In the sequel of the paper, we shall give characterize the affine
input map of the definition \ref{def:IPHS}, in such a way that the\emph{
it represents the irreversible entropy creation at the interface between
the system and its environment} in the same way as the pseudo-Poisson
bracket represents the irreversible entropy creation in the system.
Therefore we shall use a procedure that derives the port maps (the
conjugated input and output maps) by embedding the system into a IPHS
which includes both the system and its environment and then restricting
it.

\section{Port maps associated for interfaces
with irreversible phenomena}\label{sec:Irreversible-Port-maps}

In this section, we suggest to derive the port maps, that means the
input and output maps of the two conjuguated port variables, First,
we shall consider reversible explicit Port Hamiltonian Systems \cite{maschkeNOLCOS92}.
Historically, they have been derived from circuit and more generally
network-type models such as bond graphs \cite{maschkeNOLCOS92,maschkeIEEE_CAS95,maschkeJFI92}
which has been recently formalized in more general way \cite{schaft_SCL_13_PHSgraphs,Kotyczka_AT_17}.
In this paper, we shall depart from this approach and derive the port
maps from the definition of reversibel or Irrversible Hamiltonian
Systems. Therefore, we shall embed the environment as part of the
system and then restrict the total model to recover the port variables
and the Port Hamiltonian formulation. In the first step, in order
to introduce the procedure, we shall consider the reversible case
where the total model is Hamiltonian (with respect to a skew-symmetric
bracket not necessarily satisfying the Jacobi identities) and in the
second step we shall apply the procedure to the Irreversible Hamiltonian
Systems.

\subsection{Port maps for reversible Port Hamiltonian Systems }

Consider a Hamiltonian System defined on the state space which is
the product $\mathcal{X}\times\Xi$ of the space $\mathcal{X}=\mathbb{R}^{n}$
of energy variables $x\left(t\right)\in\mathcal{X}$ and the space
$\Xi$ of the environment variables $\xi\left(t\right)\in\Xi$ . Consider
$J\left(x\right)\in\mathbb{R}^{n\times n}$  a skew-symmetric matrix
(which does not satisfies in general, the Jacobi identitites), depending
on the energy variable $x$ only, and a matrix field $g\left(x\right)\in\mathbb{R}^{n\times m}$
also depending on the energy variable $x$ only. Then define the Hamiltonian
system 
\begin{equation}
\left(\begin{array}{c}
\frac{dx}{dt}\\
\frac{d\xi}{dt}
\end{array}\right)=\left(\begin{array}{cc}
J\left(x\right) & g\left(x\right)\\
-g\left(x\right)^{\top} & 0
\end{array}\right)\left(\begin{array}{c}
\frac{\partial H_{tot}}{\partial x}\\
\frac{\partial H_{tot}}{\partial\xi}
\end{array}\right)\label{eq:ExtendedHamiltonianSystem}
\end{equation}
 generated by a separable Hamiltonian $H_{tot}\left(x,\,\xi\right)=H\left(x\right)+H_{c}\left(\xi\right)$
and defined with respect to the skew-symmetric matrix 
\begin{equation}
J_{e}\left(x\right)=\left(\begin{array}{cc}
J\left(x\right) & g\left(x\right)\\
-g\left(x\right)^{\top} & 0
\end{array}\right)\label{eq:ExtendedStrctureMatrixHam}
\end{equation}
 Note that this matrix is the structure matrix of a pseudo-Poisson
bracket defined by (\ref{eq:PoissonBracket}) and is qualified as
pseudo-Poisson bracket because we do not require that it satisfies
the Jacobi identity. It is the representation, in coordinates, of
the tensor mapping the cotangent space into the tangent space (and
1-forms into vector fields) 
\begin{eqnarray*}
\Lambda:T^{*}\left(\mathcal{X}\times\Xi\right) & \rightarrow & T\left(\mathcal{X}\times\Xi\right)\\
\omega & \longmapsto & X
\end{eqnarray*}
 according to (\ref{eq:ExtendedHamiltonianSystem}).

In order to represent a control, let us consider the subspace of Hamiltonian
function $H_{c}\left(\xi\right)\in C^{\infty}\left(\mathbb{R}^{m}\right)$
is the Hamiltonian function being linear in the input: 
\begin{equation}
H_{c}\left(\xi\right)=u^{\top}\xi\label{eq:FamilyHamFunctions}
\end{equation}
 where $u\in\mathbb{R}^{m}$ are the control variables. It is clear
that the set of linear functions $L\left(\mathbb{R}^{m},\mathbb{R}\right)$
in (\ref{eq:FamilyHamFunctions}) and their differentials define a
linear vector space which is diffeomorphic to the input space $U=\mathbb{R}^{m}$.
Denote the dual to the input space $Y=U^{*}=\mathbb{R}^{m}$ identifying
$\mathbb{R}^{m}$ with its dual. Then, $\mathcal{X}\times\Xi$ may
be restricted to the space $\mathcal{X}\times L\left(\mathbb{R}^{m},\mathbb{R}\right)$
and the tensor $\Lambda$ may be identified with a linear map $\Pi:T^{*}\mathcal{X}\times U\rightarrow T\mathcal{X}\times Y$,
leading to the Port Hamiltonian System \cite{maschkeNOLCOS92} defined
on the space of the energy variables $x\in\mathcal{X}$ with port
variables $\left(u,y\right)\in U\times Y$
\begin{equation}
\left(\begin{array}{c}
\frac{dx}{dt}\\
-y
\end{array}\right)=\left(\begin{array}{cc}
J\left(x\right) & g\left(x\right)\\
-g\left(x\right)^{\top} & 0
\end{array}\right)\left(\begin{array}{c}
\frac{\partial H}{\partial x}\\
u
\end{array}\right)\label{eq:PortHamiltonianSystem}
\end{equation}

Note that this construction is very similar to the one used for deriving
Lyapunov function for forced Port Hamiltonian Systems \cite{maschkeIEEE_AC00}
and for the Control by Interconnection method \cite{OrtegaAUTOMATICA02}.

\subsection{Port maps for Irreversible Port Hamiltonian Systems}

Let s now perform the same procedure, now considering the drift dynamics
of an Irreversible Port Hamiltonian Systems (\ref{eq:IPHS}) defined
on the state space which is the product $\mathcal{X}\times\Xi$ of
the space $\mathcal{X}=\mathbb{R}^{n}$ of energy variables $x\left(t\right)\in\mathcal{X}$
and the space $\Xi$ of the environment variables $\xi\left(t\right)\in\Xi$
. And let us model solely the interaction between the energy variables
and the environment variables by considering the following \emph{anti-diagonal}
Poisson structure matrix 
\begin{equation}
J_{\textrm{port}}=\left(\begin{array}{cc}
0 & g\\
-g^{\top} & 0
\end{array}\right)\label{eq:Jport}
\end{equation}
 with \emph{constant} $g\in\mathbb{R}^{n\times m}$ . And consider
the strictly positive function $\gamma_{\textrm{port}}\left(x,\tfrac{\partial H}{\partial x},u\right)$
associated with the constitutive relation of the irreversible phenomenon
taking place at the interface between the two subsystems. As the drift
dynamics of an Irreversible Port Hamiltonian System is defined by
two functions, the total Hamiltonian and the total entropy function,
let us now consider the following functions, both linear in the environment
state variable $\xi$: the Hamiltonian function
\[
H_{\textrm{tot}}\left(x,\xi\right)=H\left(x\right)+u^{\top}\xi\;,\quad u\in\mathbb{R}^{m}
\]
 and the entropy function 
\[
S_{\textrm{tot}}\left(x,\xi\right)=S\left(x\right)+\tau^{\top}\xi\;,\quad\tau\in\mathbb{R}^{m}
\]
Note that this time, by construction two independent variables arize:
$u\in\mathbb{R}^{m}$ associated with the Hamiltonian and $\tau\in\mathbb{R}^{m}$
associated with the entropy function. The Poisson bracket, giving
the driving force of the irreversible phenomenon at the interface
is then 
\begin{eqnarray*}
\left\{ S_{\textrm{tot}},H_{\textrm{tot}}\right\} _{J_{\textrm{port}}} & = & \left(\begin{array}{c}
\frac{\partial S}{\partial x}\\
\tau
\end{array}\right)^{\top}J_{port}\left(\begin{array}{c}
\frac{\partial H}{\partial x}\\
u
\end{array}\right)\\
 & = & \left[\left(g^{\top}\frac{\partial S}{\partial x}\right)^{\top}u-\tau^{\top}\left(g^{\top}\frac{\partial H}{\partial x}\right)\right]
\end{eqnarray*}
Writing the Irreversible Port Hamiltonian drift dynamics associated
with an interface subject to an irreversible phenomenon, one obtains
hence the irreversible port maps 
\begin{eqnarray}
\frac{dx}{dt} & = & \gamma_{\textrm{port}}\left(x,\tfrac{\partial H}{\partial x},u\right)\left[\left(g^{\top}\frac{\partial S}{\partial x}\right)^{\top}u-\tau^{\top}\left(g^{\top}\frac{\partial H}{\partial x}\right)\right]g\,u\nonumber \\
\label{eq:IrrPortMap}\\
y & = & \gamma_{\textrm{port}}\left(x,\tfrac{\partial H}{\partial x},u\right)\left[\left(g^{\top}\frac{\partial S}{\partial x}\right)^{\top}u-\tau^{\top}\left(g^{\top}\frac{\partial H}{\partial x}\right)\right]g^{\top}\frac{\partial H}{\partial x}\nonumber 
\end{eqnarray}
It is immediately seen that these conjuguated port maps are the port
map of the reversible Port Hamiltonian System (\ref{eq:PortHamiltonianSystem}),
multiplied by the product 
\[
\gamma_{\textrm{port}}\left(x,\tfrac{\partial H}{\partial x},u\right)\left\{ S_{\textrm{tot}},H_{\textrm{tot}}\right\} _{J_{\textrm{port}}}
\]
 And it may be seen also that the input map might be much more general
that affine in the control variable $u$ and that the conjuguated
output map may depend on the input.

Finally, there appears an additional independent (input) variable
$\tau$ , associated with the definition of an entropy function for
the environment. It appears that for physical systems, the Hamiltonian
may be chosen as being the total energy and that the total entropy
function may not be chosen independently and may be chosen as the
sum of the entropy of each compartment \cite{Ramirez_ChemEngSci13,Ramirez_EJC13}.
In this case, the independent variable $\tau$ is fixed to $\tau=1$
and only the input $u$ remains.

\subsection{Irreversible Port Hamiltonian Systems with irreversible interface}

Using the port map (\ref{eq:IrrPortMap}), one may now define a Irreversible
Port Hamiltonian Systems with irreversible port maps as follows.
\begin{Definition}
\label{def:IPHS-IrrPortMap} An\emph{ Irreversible Port Hamiltonian
Systems with irreversible port map}, is the nonlinear control system
\begin{eqnarray}
\frac{dx}{dt} & = & \gamma\left(x,\tfrac{\partial H}{\partial x}\right)\left\{ S,H\right\} _{J}J\frac{\partial H}{\partial x}(x)\nonumber \\
 & + & \gamma_{\textrm{port}}\left(x,\tfrac{\partial H}{\partial x},u\right)\left\{ S_{\textrm{tot}},H_{\textrm{tot}}\right\} _{J_{\textrm{port}}}g\,u\label{eq:IPHS_IrrPortMap}\\
y & = & \gamma_{\textrm{port}}\left(x,\tfrac{\partial H}{\partial x},u\right)\left\{ S_{\textrm{tot}},H_{\textrm{tot}}\right\} _{J_{\textrm{port}}}g^{\top}\frac{\partial H}{\partial x}\nonumber 
\end{eqnarray}
where $x\left(t\right)\in\mathbb{R}^{n}$ is the state vector, $u\left(t\right)\in\mathbb{R}^{m}$
is the control input, and defined by 

(i) two (smooth) real functions called \emph{Hamiltonian function}
$H(x)\in C^{\infty}(\mathbb{R}^{n})$ and \emph{entropy function}
$S(x)\in C^{\infty}(\mathbb{R}^{n})$, 

(ii) the skew-symmetric \emph{structure matrix} $J\in\mathbb{R}^{n}\times\mathbb{R}^{n}$
defining a Poisson bracket (\ref{eq:PoissonBracket})

(iii) a real function $\gamma(x,\frac{\partial U}{\partial x})=\hat{\gamma}(x):C^{\infty}(\mathbb{R}^{n})$,
strictly positive function of the states and co-states

(iv) the \emph{port map} defined by the matrix $g\in\mathbb{R}^{n\times m}$
, the vector $\tau\in\mathbb{R}^{m}$ and the strictly positive function
$\gamma_{\textrm{port}}\left(x,\tfrac{\partial H}{\partial x},u\right)$
associated with constitutive relation at the interface of the system.
\end{Definition}

Let us now write the energy and entropy balance equations and therefore
recall the matrix formulation of the IPHS \begin{scriptsize}
\begin{eqnarray}
\left(\begin{array}{c}
\frac{dx}{dt}\\
-y
\end{array}\right) & = & \left[\gamma\left(x,\tfrac{\partial U}{\partial x}\right)\frac{\partial S}{\partial x}^{\top}\left(x\right)J\frac{\partial H}{\partial x}\left(x\right)\left(\begin{array}{cc}
J & 0\\
0 & 0
\end{array}\right)\right.\nonumber \\
\nonumber \\
 & + & \left.\gamma_{\textrm{port}}\left(x,\tfrac{\partial H}{\partial x},u\right)\left[\left(\begin{array}{c}
\frac{\partial S}{\partial x}\\
\tau
\end{array}\right)^{\top}J_{port}\left(\begin{array}{c}
\frac{\partial H}{\partial x}\\
u
\end{array}\right)\right]J_{port}\right]\left(\begin{array}{c}
\frac{\partial H}{\partial x}\\
u
\end{array}\right)\label{eq:IPHS_Matrix}
\end{eqnarray}
\end{scriptsize}

Multiplying (\ref{eq:IPHS_Matrix}) from the left by $\left(\begin{array}{cc}
\frac{\partial H}{\partial x}^{\top} & u^{\top}\end{array}\right)$ , by the skew symmetry of the structure matrices, one obtains the
energy balance equation
\[
\frac{dH}{dt}-y^{\top}u=0
\]
 which corresponds to a lossless system with be impedance type of
supply function. 

Multiplying (\ref{eq:IPHS_Matrix}) from the left by $\left(\begin{array}{cc}
\frac{\partial S}{\partial x}^{\top} & \tau^{\top}\end{array}\right)$ one obtains the entropy balance equation 
\begin{eqnarray*}
\frac{dS}{dt}-\tau^{\top}y & = & \underbrace{\gamma\left(x,\tfrac{\partial U}{\partial x}\right)\left\{ S,U\right\} _{J}^{2}}_{=\sigma_{int}\geq0}\\
 & + & \underbrace{\gamma_{\textrm{port}}\left(x,\tfrac{\partial H}{\partial x},u\right)\left\{ S_{\textrm{tot}},H_{\textrm{tot}}\right\} _{J_{\textrm{port}}}^{2}}_{=\sigma_{port}\geq0}
\end{eqnarray*}
where the right-hand side is positive and corresponds to the irreversible
entropy creation due to the irreversible processes in the system $\sigma_{int}$
and at its port $\sigma_{port}$ and the term $\tau^{\top}y$ corresponds
to the entropy flowing out the environment (to the system).

\subsection{Heat conduction in 2 compartement system}

Consider again the example of the two compartements system considered
in the preceeding section. As the compartment 2 is in interaction
with the external thermostat through a heat conducting wall, as the
two compartement are, the structure matrix associated with this wall
is the symplectic matrix $J_{port}=\left[\begin{smallmatrix}0 & 1\\
-1 & 0
\end{smallmatrix}\right]$ and the postive function associated with the heat conduction relation
is $\gamma_{port}\left(x\right)=\frac{\lambda}{T_{2}u}$ where the
input is the thermostat temperature $u=T_{e}$ and fixing $\tau=1$.

The Irreversible Port Hamiltonian Systems with irreversible port map
is then composed of the state equation 
\begin{eqnarray*}
\frac{d}{dt}\begin{bmatrix}x_{1}\\
x_{2}
\end{bmatrix} & = & \frac{\lambda}{T_{1}T_{2}}\left(T_{2}-T_{1}\right)\begin{bmatrix}0 & 1\\
-1 & 0
\end{bmatrix}\begin{bmatrix}\tfrac{\partial H}{\partial x_{1}}\\
\tfrac{\partial H}{\partial x_{2}}
\end{bmatrix}\\
 & + & \frac{\lambda_{e}}{uT_{2}}\left(u-T_{2}\right)\begin{bmatrix}0\\
1
\end{bmatrix}u
\end{eqnarray*}

where the input map may be simplified to become affine in the input
\[
\frac{\lambda_{e}}{uT_{2}}\left(u-T_{2}\right)\begin{bmatrix}0\\
1
\end{bmatrix}u=\frac{\lambda_{e}}{T_{2}}\left(u-T_{2}\right)\begin{bmatrix}0\\
1
\end{bmatrix}
\]
 and the output equation 
\[
y=\frac{\lambda_{e}}{uT_{2}}\left[\left(u-T_{2}\right)\right]T_{2}=\frac{\lambda_{e}}{u}\left[\left(u-T_{2}\right)\right]
\]
 which is the entropy flux leaving the environment ( \emph{with respect
to the temperature }$u=T_{e}>0$\emph{ of the environment}).

The entropy balance equation is then 
\begin{eqnarray*}
\frac{dS}{dt} & = & \frac{\lambda}{T_{1}T_{2}}\left(T_{2}-T_{1}\right)^{2}\\
 & + & \frac{\lambda_{e}}{uT_{2}}\left(u-T_{2}\right)^{2}+y
\end{eqnarray*}
 and the definition of the input then naturally makes appear the irreversible
entropy creation at the interface associated with the port of the
system as in (\ref{eq:IrrevesibelEntrCrea_WallExt}).

\section{Conclusion}

In this paper, we have suggested a novel definition of the port maps
of Irreversible Port Hamiltonian Systems and have derived it by embedding
the environment into a composed Irreversible Port Hamiltonian system
and then restricting the definition of the energy and entropy function
associated with the environment to linear functions. As a result,
we have defined a novel port map consisting in a nonlinear input map
which is defined by the constitutive relation of an irreversible phenomenon
and have derived a conjuguated output which allows to recover a entropy
balance equation including not only the irreversible entropy creation
at the interface (the port) of the system in addition to the entropy
creation etrm due to internal irreversible phenomena. In this way,
we have generalized and given a physical meaningful structure to the
input map affine in the input suggested in \cite{Ramirez_ChemEngSci13}.
This construction has been illustrated on the very elementary example
of the system of enetropy balance equations of a 2-compartement system
subject to heat conduction and in interaction with a thermostat.

As a result of associating not only a (linear) Hamiltonian but also
an entropy function to the environment, in addition to the port input
variable an additional independent parameter has been introduced which
could eventually be also viewed as an input. Hence the obtained input-map
is more a linear map as for reversible Port Hamiltonian systems, but
non-linear in two independent variables (input port variables). In
future work, we shall elucidate this structure with a more geometric
approach by extending the formulation of the drift term of Irreversible
Port Hamiltonian systems presented in a companion paper at this conference
\cite{Kirchhoff_submIFAC_WC22}.

\bibliographystyle{alpha}

\end{document}